\documentclass[12pt]{amsart}

\usepackage{amsmath}
\usepackage{amssymb}
\usepackage{array}

\renewcommand{\atop}[2]{%
\genfrac{}{}{0pt}{}{#1}{#2}}

\newtheorem{theorem}{Theorem}
\newtheorem{lemma}{Lemma}
\newtheorem{corollary}{Corollary}

\theoremstyle{definition}

\renewcommand{\atop}[2]{%
\genfrac{}{}{0pt}{}{#1}{#2}}

\begin{document}

\title{Rational approximations for the quotient of \\ gamma values}

\author{Kh.~Hessami Pilehrood$^1$}

\address{\begin{flushleft} Mathematics Department, Faculty of Basic
Sciences, Shahrekord University, Shahrekord, P.O. Box 115, Iran. \end{flushleft}}
\thanks{$^1$ This research was in part supported by a grant
from IPM (No. 89110024)}

\author{T.~Hessami Pilehrood$^2$}
\address{\begin{flushleft} School of Mathematics, Institute for Research  in Fundamental Sciences
(IPM), P.O.Box 19395-5746, Tehran, Iran \end{flushleft}}


\email{hessamik@ipm.ir, hessamit@ipm.ir, hessamit@gmail.com}
\thanks{$^2$ This research was in part supported by a grant
from IPM (No. 89110025)}


\date{}


\begin{abstract}
In this paper, we continue to study properties of rational approximations to Euler's constant
and values of the Gamma function defined by linear recurrences, which were found recently by
A.~I.~Aptekarev and T.~Rivoal. Using multiple Jacobi-Laguerre orthogonal polynomials we present rational
approximations to the  quotient of values of the Gamma function at rational points.
As a limit case of our result, we obtain new explicit formulas for numerators and denominators of the Aptekarev
approximants to Euler's constant.
\end{abstract}

\maketitle

\section{Introduction}

Recently, A.~I.~Aptekarev and his collaborators in a series of papers  \cite{ap,tu} proposed a sequence of
rational approximations $p_n/q_n$ converging to Euler's constant
$$
\gamma=\lim_{k\to\infty}\left(1+\frac{1}{2}+\cdots+\frac{1}{k}-\ln k\right).
$$
sub-exponentially.
The numerators $p_n$ and the denominators $q_n$ of the approximations are positive integers
generated by the following recurrence relation:
\begin{equation}
\begin{split}
(16n-15)q_{n+1}&=(128n^3+40n^2-82n-45)q_n \\
&-n^2(256n^3-240n^2+64n-7)q_{n-1}
+n^2(n-1)^2(16n+1)q_{n-2}
\label{eq00}
\end{split}
\end{equation}
with initial conditions
\begin{equation}
\begin{array}{ccc}
p_0=0, \qquad & \qquad p_1=2, \qquad & \qquad p_2=31, \\
q_0=1,  \qquad &  \qquad q_1=3, \qquad &  \qquad q_2=50
\label{ini}
\end{array}
\end{equation}
and having the following asymptotics:
\begin{align}
q_n&=(2n)!\frac{e^{\sqrt{2n}}}{\sqrt[4]{n}}\left(\frac{1}{\sqrt{\pi}(4e)^{3/8}}+O(n^{-1/2})\right),\nonumber\\[3pt]
p_n-\gamma q_n&=(2n)!\frac{e^{-\sqrt{2n}}}{\sqrt[4]{n}}\left(\frac{2\sqrt{\pi}}{(4e)^{3/8}}+O(n^{-1/2})\right)
\label{eq01}
\end{align}
In \cite{he}, explicit formulas for the sequences $p_n$ and $q_n$ were given
in terms of hypergeometric sums
\begin{equation}
q_n=\sum_{k=0}^n\binom{n}{k}^2(n+k)!,
\label{znam}
\end{equation}
$$
p_n=\sum_{k=0}^n\binom{n}{k}^2(n+k)!H_{n+k}
-2\sum_{k=1}^n\sum_{m=0}^k\sum_{l=0}^{n-k}\frac{(-1)^{m+k}}{k}\binom{n}{k+l}\binom{k}{m}\binom{n-k}{l}(m+n+l)!,
$$

\vspace{0.2cm}

\noindent where $H_n$ is an $n$-th harmonic number defined by
$$
H_n(\alpha)=\sum_{k=1}^n\frac{1}{k+\alpha}, \quad  n\ge 1, \quad
H_0(\alpha)=0, \quad H_n=H_n(0).
$$
These representations provide a more profound insight on arithmetical properties of these sequences
that can not be deduced from the recurrence equation (\ref{eq00}).
From (\ref{znam}) it follows that
$q_n$ and $p_n$ are integers divisible by $n!$ and $n!/D_n$ respectively,
where  $D_n$ is the least common multiple of the numbers $1,2,\ldots,n.$
This implies that the coefficients of the linear form
(\ref{eq01})
can be canceled out  by the big common factor $\frac{n!}{D_n},$ but it is still not sufficient  to prove
the irrationality of $\gamma.$
The remainder of the above approximation is given by the integral \cite{at}
$$
\int_0^{\infty}Q_n^{(0,0,1)}(x)\,e^{-x}\ln(x)\,dx=p_n-\gamma q_n,
$$
where for $a_1, a_2>-1, b>0,$
\begin{equation}
Q_n^{(a_1,a_2,b)}(x)=
\frac{1}{n!^2}\prod_{j=1}^2\bigl[w_j^{-1}\frac{d^n}{dx^n}w_jx^n\bigr](1-x)^{2n}
\label{eq02}
\end{equation}
is a multiple Jacobi-Laguerre orthogonal polynomial
on $\Delta_1=[0,1]$
and $\Delta_2=[1,+\infty)$
with respect to the two weight functions
$$
w_1(x)=x^{a_1}(1-x)e^{-bx}, \qquad w_2(x)=x^{a_2}(1-x)e^{-bx}
$$
 satisfying the following orthogonality conditions \cite{as}:
\begin{equation}
\int_{\Delta_j}x^{\nu}Q_n^{(a_1,a_2,b)}(x)w_k(x)\, dx=0, \quad j,k=1,2, \,\,
\nu=0,1,\ldots,n-1.
\label{eq03}
\end{equation}
Note that, as the values of the parameters $a_1, a_2$ tend to $0$ and $b=1,$  the limit values
of $Q_n^{(0,0,1)}(x)$ satisfy the orthogonality relations (\ref{eq03}) with the weights
$w_1(x)=(1-x)e^{-x},$ $w_2(x)=(1-x)\ln(x)e^{-x}.$

In \cite{he}, the authors showed that for $a, b\in {\mathbb Q},$ $a>-1,$ $b>0$ the integrals
\begin{equation}
\frac{b^{2n+a+1}}{\Gamma(a+1)}\int_0^{\infty}Q_n^{(a,a,b)}(x)x^a\ln(x)e^{-bx}\,dx=
p_n(a,b)-(\ln(b)-\psi(a+1))q_n(a,b),
\label{pn}
\end{equation}
where $\psi(z)=\Gamma'(z)/\Gamma(z)$
defines the logarithmic derivative of the gamma function,
provide good rational approximations for $\ln(b)-\psi(a+1),$  which converge sub-exponentially
to this number
\begin{equation}
\frac{p_n(a,b)}{q_n(a,b)}-(\ln(b)-\psi(a+1))=2\pi e^{-2\sqrt{2bn}}(1+O(n^{-1/2})) \quad
\mbox{as} \,\,\,\,  n\to\infty,
\label{eq03.5}
\end{equation}
and gave explicit formulas for the denominators and numerators of the approximations
\begin{equation}
\begin{split}
q_n(a,b)&=\sum_{k=0}^n\binom{n}{k}^2(a+1)_{n+k}\,b^{n-k}, \\
\label{eq04}
p_n(a,b)&=\sum_{k=0}^n\binom{n}{k}^2(a+1)_{n+k}\,b^{n-k}H_{n+k}(a) \\
&-2\sum_{k=1}^n\sum_{m=0}^k\sum_{l=0}^{n-k}\frac{(-1)^{m+k}}{k}\binom{n}{k+l}\binom{k}{m}
\binom{n-k}{l}(a+1)_{m+n+l}\,b^{n-m-l},
\end{split}
\end{equation}
here $(a)_n$  is the Pochhammer symbol defined by $(a)_0=1$ and $(a)_n=a(a+1)\cdots (a+n-1),$
$n\ge 1.$

Recently, T.~Rivoal \cite{ri} found another example of rational approximations to the number
$\gamma+\ln(x),$ $x\in {\mathbb Q},$ $x>0,$ by using multiple Laguerre polynomials
$$
A_n(x)=\frac{e^x}{n!^2}(x^n(x^ne^{-x})^{(n)})^{(n)}=e^x\cdot
{}\sb 2F\sb
 {2}\left(\left.\atop{n+1,
n+1}{1, \quad 1}\right|-x\right)
$$
(for the validity of the last equality see \cite[\S 4.3]{lu}).
Here ${}\sb pF\sb {q}$ is the generalized hypergeometric function defined by the series
$$
{}\sb pF\sb
 {q}\left(\left.\atop{\alpha_1, \ldots,
\alpha_p}{\beta_1, \ldots \beta_q}\right|z\right)=\sum_{k=0}^{\infty}
\frac{(\alpha_1)_k\cdots (\alpha_p)_k}{(\beta_1)_k\cdots (\beta_q)_k}
\frac{z^k}{k!},
$$
which converges for all finite $z$ if $p\le q,$ converges for $|z|<1$ if $p=q+1,$
and diverges for all $z,$ $z\ne 0$ if $p>q+1.$ The series ${}\sb pF\sb {q}$ terminates
and, therefore, is a polynomial in $z$ if a numerator parameter is a negative integer or zero,
provided that no denominator parameter is a negative integer or zero.

The polynomial $A_n(x)$ as well as the hypergeometric function ${}\sb 2F\sb
 {2}\bigl(\bigl.\atop{n+1,
n+1}{1, \quad 1}\bigr|-x\bigr)$ satisfies a certain third-order linear recurrence relation with
polynomial coefficients in $n$ and $x.$ The polynomials $A_n(x)$ are orthogonal with respect
to the weight functions $\omega_1(x)=e^{-x},$ $\omega_2(x)=e^{-x}\ln(x)$ on $[0,+\infty)$, so that
\begin{equation}
\int_0^{\infty}x^{\nu}A_n(x)\omega_k(x)\,dx=0, \qquad k=1,2, \,\,  \nu=0,1,\ldots, n-1.
\label{eq05}
\end{equation}
This implies that $A_n(x)$ is a common denominator in type II  Hermite-Pad\'e approximation for the two
Stieltjes functions
$$
{\mathcal E}_1(x)=\int_0^{\infty}\frac{e^{-t}}{x-t}\,dt, \qquad
{\mathcal E}_2(x)=\int_0^{\infty}\frac{e^{-t}\ln(t)}{x-t}\,dt
$$
at $x=\infty.$ The remainders of this approximation are given by the integrals
\begin{equation*}
\begin{split}
R_{1,n}(x):=\int_0^{\infty}\frac{A_n(t)}{x-t}e^{-t}\,dt&=A_n(x){\mathcal E}_1(x)+B_n(x), \\
R_{2,n}(x):=\int_0^{\infty}\frac{A_n(t)}{x-t}\ln(t)e^{-t}\,dt&=A_n(x){\mathcal E}_2(x)+C_n(x),
\end{split}
\end{equation*}
which in view of the orthogonality relations (\ref{eq05}) satisfy the same recurrence as the polynomial
$A_n(x)$ (see \cite[\S 5.2]{ri}). By linearity, it is easily seen that the remaining polynomials
$B_n(x), C_n(x),$ and
the function
$$
R_n(x)=\ln(x) R_{1,n}(x)-R_{2,n}(x)=\int_0^{\infty}A_n(t)\frac{\ln(x)-\ln(t)}{x-t}e^{-t}\, dt
$$
are also solutions of the same recurrence. Suppose that this third-order recurrence has the form
$$
U_{n+3}=s_nU_{n+2}+t_nU_{n+1}+r_nU_n.
$$
Then rational approximations to $\ln(x)+\gamma$ are given \cite{ri} by a so called ``determinantal''
sequence
\begin{equation}
{\mathcal S}_n(x)=\begin{vmatrix}
A_n(x) & R_n(x) \\
A_{n+1}(x) & R_{n+1}(x)
\end{vmatrix}
=Q_n(x)(\ln(x)+\gamma)-P_n(x), \quad Q_n(x), P_n(x)\in{\mathbb Q}[x],
\label{10.1}
\end{equation}
which along with the sequences $Q_n(x)$ and $P_n(x)$ satisfies a ``determinantal'' recurrence of the form
$$
V_{n+3}=-t_{n+1}V_{n+2}-s_nr_{n+1}V_{n+1}+r_nr_{n+1}V_n.
$$
The sequence $P_n(x)/Q_n(x)$ converges to $\ln(x)+\gamma$ sub-exponentially
$$
\left|\ln(x)+\gamma-\frac{P_n(x)}{Q_n(x)}\right|\le c_0(x) e^{-9/2x^{1/3}n^{2/3}+3/2x^{2/3}n^{1/3}},
\quad Q_n(x)=O(e^{3x^{1/3}n^{2/3}-x^{2/3}n^{1/3}}),
$$
and for  denominators of  coefficients of the polynomials $P_n$ and $Q_n$ one has
$$
n!\,Q_n(x)\in{\mathbb Z}[x], \qquad D_n n!\,P_n(x)\in{\mathbb Z}[x]
$$
(the first inclusion was formulated in \cite{ri} as a denominator conjecture, and both
inclusions were proved in \cite{he}).

In \cite{ri1}, the similar scheme was applied to the Laguerre polynomials
$$
A_{n,a}(x)=\frac{e^x}{n!^2}(x^{n-a}(x^{n+a}e^{-x})^{(n)})^{(n)}=
\frac{(a+1)_n}{n!}\,e^x
{}\sb 2F\sb
 {2}\left(\left.\atop{a+n+1,
n+1}{a+1, \, 1}\right|-x\right),
$$
which are orthogonal on $[0,+\infty)$ with respect to the weights $e^{-x}$ and $x^ae^{-x},$
to get rational approximations for $\Gamma(a+1)/x^a,$ where $a,x\in {\mathbb Q}, a>-1, x>0.$
In this case, rational approximations are defined by the ``determinantal'' sequence
\begin{equation}
{\mathcal S}_n(a,x)=\begin{vmatrix}
A_{n,a}(x) & R_{n,a}(x) \\
A_{n+1,a}(x) & R_{n+1,a}(x)
\end{vmatrix}
=Q_n(a,x)\,\Gamma(a+1)-P_n(a,x)\,x^a,
\label{10.5}
\end{equation}
where
$$
R_{n,a}(x)=x^aR_{n,a,0}(x)-R_{n,a,a}(x) \quad\mbox{and}\quad
R_{n,a,c}(x)=\int_0^{\infty}\frac{A_{n,a}(t)}{x-t}\,t^ce^{-t}\,dt.
$$
The sequences $P_n(a,x), Q_n(a,x), {\mathcal S}_n(a,x)$ are solutions of the third-order
recurrence
$$
C_3(n,a,x)V_{n+3}+C_2(n,a,x)V_{n+2}+C_1(n,a,x)V_{n+1}+C_0(n,a,x)V_n=0,
$$
where the coefficients $C_j(n,a,x),$ $j=0,1,2,3$  are polynomials with integer coefficients
in $n, a, x$ of degree 16 in $n.$ Moreover, one has \cite{ri1}
\begin{equation*}
\begin{split}
|Q_n(a,x) \Gamma(a+1)-P_n(a,x) x^a|&\le
c_1(a,x)\, e^{-3/2x^{1/3}n^{2/3}+1/2x^{2/3}n^{1/3}}n^{2a/3-2}, \\
Q_n(a,x)&=O(e^{3x^{1/3}n^{2/3}-x^{2/3}n^{1/3}}n^{2a/3-2})
\end{split}
\end{equation*}
with the following inclusions proved in \cite{he}:
$$
(n+1)^2\cdot n!\,\mu_a^{2n+1}({\rm den}\,x)^{n+1}P_n(a,x), \,\,
(n+1)^2\cdot n!\,\mu_a^{3n+1}({\rm den}\,x)^{n+1}Q_n(a,x)\in {\mathbb Z},
$$
where
\begin{equation}
\mu_a^n=({\rm den}\,a)^n\cdot\prod_{p|{\rm den}\,a}p^{[\frac{n}{p-1}]}
\label{eq06}
\end{equation}
and ${\rm den}\,a\in {\mathbb N}$ is the denominator of a simplified fraction $a.$

In this paper, we continue to study  properties of the above constructions
and get rational approximations for the quotient of values
of the Gamma function at rational points.
The main result of this paper is as follows.
\begin{theorem} \label{t1}
Let $a_1, a_2, b\in {\mathbb Q},$ $a_1-a_2\not\in {\mathbb Z},$ $a_1, a_2>-1, b>0.$
For $n=0,1,2,\ldots,$ define a sequence of rational numbers
\begin{equation*}
q_n(a_1,a_2,b)=\sum_{k=0}^n\binom{n+a_1-a_2}{k}\binom{n+a_2-a_1}{n-k}(a_2+1)_{n+k}b^{n-k},
\end{equation*}
Then
$$
\mu_{a_2-a_1}^n\mu_{a_2}^{2n}\,q_n(a_1,a_2,b) 
\in n!\,{\mathbb Z}[b],
$$

\vspace{0.4cm}

\noindent where $\mu_a^n$ is defined in {\rm (\ref{eq06})}, and the following asymptotic formulas are valid:
\begin{equation}
q_n(a_1,a_2,b)\frac{\Gamma(a_2+1)/b^{a_2}}{\Gamma(a_1+1)/b^{a_1}} - q_n(a_2,a_1,b)=(2n)!
\frac{e^{-\sqrt{2bn}}}{n^{1/4-(a_1+a_2)/2}}(c_1+O(n^{-1/2})),
\label{eq06.5}
\end{equation}
$$
q_n(a_1,a_2,b)=(2n)!\frac{e^{\sqrt{2bn}}}{n^{1/4-(a_1+a_2)/2}}(c_2+O(n^{-1/2})) \quad
\mbox{as}\,\,\, n\to\infty
$$
with
\begin{equation}
c_1=\frac{2\sin(\pi(a_2-a_1))\Gamma(a_2+1)}{b^{a_2-a_1}\Gamma(a_1+1)}\,c_2,
\quad
c_2=\frac{2^{(a_1+a_2)/2-3/4}}{b^{1/4+(a_1-a_2)/2}e^{3b/8}\sqrt{\pi}\Gamma(a_2+1)}.
\label{eq07}
\end{equation}
\end{theorem}
The linear forms (\ref{eq06.5}) do not allow one to prove the irrationality of the quotient
$\Gamma(a_2+1)/(\Gamma(a_1+1)b^{a_2-a_1}),$ however they present good rational approximations
to this number
$$
\frac{\Gamma(a_2+1)/b^{a_2}}{\Gamma(a_1+1)/b^{a_1}} - \frac{q_n(a_2,a_1,b)}{q_n(a_1,a_2,b)}=
\frac{c_1}{c_2}\,\,e^{-2\sqrt{2bn}}\,
\bigl(1+O(n^{-1/2})\bigr).
$$
For an overview of known results on arithmetical nature of Gamma values, see \cite{ri1}
and the references given there.

Application of Zeilberger's algorithm of creative telescoping \cite[Ch.~6]{pwz} shows
that both sequences $q_n(a_1,a_2,b)$ and $q_n(a_2,a_1,b)$ satisfy a third-order linear recurrence
with polynomial coefficients in $n, a_1, a_2, b$ of degree not exceeding 10 in $n.$
We do not write down this recurrence equation for arbitrary $a_1, a_2, b,$ since it is quite
complicated in its general form. We indicate  this recurrence only in several particular cases.

Setting $b=1,$ $a_2=-1/2,$ $a_1=-2/3$  and then $a_1=-3/4$ we deduce the following results.
\begin{corollary}
The sequences
\begin{equation*}
p_n=\sum_{k=0}^n\binom{n+\frac{1}{6}}{k}\binom{n-\frac{1}{6}}{n-k}\Bigl(\frac{1}{3}\Bigr)_{n+k}
\quad\mbox{and}\quad\,\,
q_n=\sum_{k=0}^n\binom{n-\frac{1}{6}}{k}\binom{n+\frac{1}{6}}{n-k}\Bigl(\frac{1}{2}\Bigr)_{n+k}
\end{equation*}
satisfy the third-order linear recurrence
\begin{equation*}
\begin{split}
&7776(n+3)(n+2)(1728n^4+6228n^3+7986n^2+4538n+963)f_{n+3}\\
&=
36(2985984n^7+33654528n^6+155405952n^5+379069848n^4+525871008n^3 \\
&+415828470n^2+175094660n+30377295)(n+2)f_{n+2} \\
&-6(6n+11)(6n+13)(995328n^8+11384064n^7+55911168n^6+154077840n^5 \\
&+260388504n^4+275942444n^3
+178728654n^2+64534088n+9914615)f_{n+1} \\
&+(6n+11)(6n+5)(2n+1)(3n+1)(6n+13)(6n+7)(1728n^4+13140n^3 \\
&+37038n^2+46106n+21443)f_n
\end{split}
\end{equation*}
with the initial values
$p_0=1,$ $p_1=43/54,$ $p_2=12871/1458$ and $q_0=1,$ $q_1=29/24,$ $q_2=5149/384,$
and provide rational approximations for the following quotient:
$$
\frac{\Gamma(1/2)}{\Gamma(1/3)}-\frac{p_n}{q_n}=e^{-2\sqrt{2n}}\Bigl(\frac{\sqrt{\pi}}{\Gamma(1/3)}+O(n^{-1/2})\Bigr)
\qquad\mbox{as}\quad n\to\infty.
$$
\end{corollary}
\begin{corollary}
The sequences
\begin{equation*}
p_n=\sum_{k=0}^n\binom{n+\frac{1}{4}}{k}\binom{n-\frac{1}{4}}{n-k}\Bigl(\frac{1}{4}\Bigr)_{n+k}
\quad\mbox{and}\quad\,\,
q_n=\sum_{k=0}^n\binom{n-\frac{1}{4}}{k}\binom{n+\frac{1}{4}}{n-k}\Bigl(\frac{1}{2}\Bigr)_{n+k}
\end{equation*}
satisfy the third-order linear  recurrence
\begin{equation*}
\begin{split}
&2048(n+2)(n+1)(2048n^4-1024n^3-304n^2+416n-117)f_{n+2} \\
&=32(4n+5)(n+1)(262144n^6+753664n^5+331776n^4-365824n^3+17168n^2 \\
&+76032n
-32229)f_{n+1}
-2(4n+3)(4n+5)(4n+1)(524288n^7+1572864n^6 \\
&+1486848n^5+389120n^4-245824n^3-117296n^2-604n-537)f_n \\
&+(4n+3)(4n-1)(2n-1)(4n+5)(4n+1)(4n-3)(2048n^4+7168n^3 \\
&+8912n^2+4928n+1019)f_{n-1}
\end{split}
\end{equation*}
with the initial values $p_0=1,$ $p_1=37/64,$ $p_2=50685/8192$ and $q_0=1,$ $q_1=19/16,$ $q_2=6525/512,$
and provide rational approximations for the following quotient:
$$
\frac{\Gamma(1/2)}{\Gamma(1/4)}-\frac{p_n}{q_n}=e^{-2\sqrt{2n}}\Bigl(\frac{\sqrt{2\pi}}{\Gamma(1/4)}+O(n^{-1/2})\Bigr)
\qquad\mbox{as}\quad n\to\infty.
$$
\end{corollary}
Setting $a_1=0,$ $a_2=a,$ $b=1$ in Theorem \ref{t1} we get
rational approximations for values of the Gamma function.
\begin{corollary}
Let  $a\in {\mathbb Q}\setminus{\mathbb Z},$ $a>-1.$ Then the sequences
\begin{equation*}
p_n=\sum_{k=0}^n\binom{n+a}{k}\binom{n-a}{n-k}(n+k)!
\quad\mbox{and}\quad\,\,
q_n=\sum_{k=0}^n\binom{n-a}{k}\binom{n+a}{n-k}(a+1)_{n+k}
\end{equation*}
satisfy the third-order recurrence
$$
\lambda_2(n,a)f_{n+2}+\lambda_1(n,a)f_{n+1}+\lambda_0(n,a)f_n+\lambda_{-1}(n,a)f_{n-1}=0,
$$
where
\begin{equation*}
\begin{split}
\lambda_2(n,a)&=(n+2)(16n^3+20n^2a+10na^2+2a^3+17n^2+14na+3a^2+n-a), \\
\lambda_1(n,a)&=-128n^6-224n^5a-112n^4a^2+20n^3a^3+42n^2a^4+16na^5+2a^6-808n^5 \\ &-1228n^4a-601n^3a^2-45n^2a^3+47na^4+11a^5
-1910n^4-2446n^3a \\
&-1085n^2a^2-176na^3-3a^4-2035n^3-2063n^2a-715na^2-91a^3-887n^2 \\
&-592na-113a^2-82n+14a, \\
\lambda_0(n,a)&=(n-a+1)(n+a+1)(256n^6+576n^5a+544n^4a^2+272n^3a^3+72n^2a^4 \\
&+8na^5+1296n^5+2384n^4a+1724n^3a^2+596n^2a^3
+92na^4+4a^5+2448n^4 \\
&+3518n^3a+1840n^2a^2+404na^3+30a^4+2185n^3+2333n^2a+817na^2+93a^3 \\
&+923n^2+668na+125a^2+146n+58a), \\
\lambda_{-1}(n,a)&=
-(n-a)(n-a+1)(n+a+1)(16n^3+20n^2a+10na^2+2a^3+65n^2 \\
&+54na+13a^2+83n+33a+34)(n+a)^2,
\end{split}
\end{equation*}
with the initial conditions
\begin{equation*}
\begin{array}{ccc}
p_0=1, \quad &  p_1=3+a, \quad &  p_2=50+33a+7a^2, \\
q_0=1,  \quad &   q_1=(a+1)(3-a^2), \quad &   q_2=(a+1)(a+2)(a^4+2a^3-12a^2-11a+50)/2,
\end{array}
\end{equation*}
and provide rational approximations for the following value of the Gamma function:
$$
\Gamma(a+1)-\frac{p_n}{q_n}=
e^{-2\sqrt{2n}}\Bigl(\frac{2\pi a}{\Gamma(1-a)}+O(n^{-1/2})\Bigr)
\qquad\mbox{as}\quad n\to\infty.
$$
\end{corollary}

Another positive feature of our considerations is the fact that we can prove a simpler
formula for the numerator $p_n(a,b)$ of rational approximations (\ref{eq03.5})
than a representation given by (\ref{eq04}).
\begin{theorem} \label{t2}
For the sequence $p_n(a,b),$ $n=0,1,2,\ldots,$ defined by {\rm (\ref{eq04})},
a more compact formula is valid
$$
p_n(a,b)=\sum_{k=0}^n\binom{n}{k}^2(a+1)_{n+k}\,b^{n-k}(H_{n+k}(a)+2H_{n-k}-2H_k).
$$
\end{theorem}
As a consequence, we get a simple  formula for the numerator $p_n$ of rational
approximations (\ref{eq01})
to Euler's constant.
\begin{corollary} \label{c1}
The sequences $p_n$ and $q_n,$ $n=0,1,2,\ldots,$ defined by the recurrence {\rm (\ref{eq00})} with initial
conditions {\rm (\ref{ini})} are given by the formulas
\begin{equation*}
\begin{split}
q_n&=\sum_{k=0}^n\binom{n}{k}^2(n+k)!, \\
p_n&=\sum_{k=0}^n\binom{n}{k}^2(n+k)!(H_{n+k}+2H_{n-k}-2H_k).
\end{split}
\end{equation*}
\end{corollary}

\vspace{0.1cm}

\section{Analytical construction}

\vspace{0.2cm}

Let $a_1, a_2, b$ be arbitrary rational numbers such that $a_1-a_2\not\in {\mathbb Z},$ $a_1, a_2>-1,$
$b>0.$ The construction of rational approximations to the numbers $\Gamma(a_2+1)/(\Gamma(a_1+1)b^{a_2-a_1})$
is based on the integral
\begin{equation}
R_n:=R_n(a_1,a_2,b)=b^{2n+1}\int_0^{\infty}
(x^{a_2}-x^{a_1})Q_n^{(a_1,a_2,b)}(x)e^{-bx}\,dx,
\label{eq12}
\end{equation}
where the polynomial $Q_n^{(a_1,a_2,b)}(x)$ is defined in (\ref{eq02}).
\begin{lemma} \label{l1}
Let $a_1, a_2, b\in {\mathbb Q},$ $a_1-a_2\not\in {\mathbb Z},$ $a_1, a_2>-1,$
$b>0.$ Then we have

$$
R_n(a_1,a_2,b)=q_n(a_1,a_2,b)\frac{\Gamma(a_2+1)}{b^{a_2}}
-q_n(a_2,a_1,b)\frac{\Gamma(a_1+1)}{b^{a_1}},
$$

\vspace{0.3cm}

\noindent where $q_n(a_1,a_2,b)$ is defined in Theorem {\rm \ref{t1}},
\begin{equation}
q_n(a_1,a_2,b)=\frac{b^{2n+a_2+1}}{\Gamma(a_2+1)}\int_0^{\infty}x^{a_2}Q_n^{(a_1,a_2,b)}
e^{-bx}\,dx
\label{eq13}
\end{equation}
and $\mu_{a_2-a_1}^n\mu_{a_2}^{2n}\,q_n(a_1,a_2,b)
\in n!\,{\mathbb Z}[b].$
\end{lemma}
{\bf Proof.} Substituting (\ref{eq02}) into (\ref{eq12}) we get
\begin{equation}
\begin{split}
R_n(a_1,a_2,b)&=\frac{b^{2n+1}}{n!^2}\int_0^{\infty}
\frac{(x^{n+a_2-a_1}(x^{n+a_1}(1-x)^{2n+1}e^{-bx})^{(n)})^{(n)}}{1-x}\,dx \\[4pt]
&-\frac{b^{2n+1}}{n!^2}\int_0^{\infty}
\frac{x^{a_1-a_2}}{1-x}(x^{n+a_2-a_1}(x^{n+a_1}(1-x)^{2n+1}e^{-bx})^{(n)})^{(n)}\,dx.
\end{split}
\label{eq14}
\end{equation}
Denoting the first integral on the right-hand side of (\ref{eq14}) by $I_1$ and applying
$n$-integrations by parts, we have
\begin{equation*}
\begin{split}
I_1&=\frac{(-1)^n b^{2n+1}}{n!}\int_0^{\infty}
\frac{x^{n+a_2-a_1}}{(1-x)^{n+1}}\bigl(x^{n+a_1}(1-x)^{2n+1}e^{-bx}\bigr)^{(n)}\,dx \\[2pt]
&=\frac{b^{2n+1}}{n!}\int_0^{\infty}
\left(\frac{x^{n+a_2-a_1}}{(1-x)^{n+1}}\right)^{(n)}
x^{n+a_1}(1-x)^{2n+1}e^{-bx}\,dx.
\end{split}
\end{equation*}
Now using the following relation for hypergeometric functions (see \cite[\S 3.4, (11)]{lu}):
\begin{equation}
\frac{d^n}{dx^n}[x^{\alpha}(1-x)^{\beta}]=(\alpha-n+1)_n x^{\alpha-n} (1-x)^{\beta-n}
{}\sb 2F\sb
 {1}\!\left(\!\left.\atop{-n,
\alpha+1+\beta-n}{\alpha+1-n}\right|x\right), \,\, n-\alpha\not\in {\mathbb N},
\label{id}
\end{equation}
with $\alpha=n+a_2-a_1,$ $\beta=-n-1,$ we get
$$
\left(\frac{x^{n+a_2-a_1}}{(1-x)^{n+1}}\right)^{(n)}=(a_2-a_1+1)_n x^{a_2-a_1}(1-x)^{-2n-1}
{}\sb 2F\sb
 {1}\!\left(\!\left.\atop{-n,
a_2-a_1-n}{a_2-a_1+1}\right|x\right)
$$
and therefore,
\begin{equation}
\begin{split}
I_1&=\frac{b^{2n+1}(1+a_2-a_1)_n}{n!}\int_0^{\infty}
x^{n+a_2}{}\sb 2F\sb
 {1}\!\left(\!\left.\atop{-n,
a_2-a_1-n}{a_2-a_1+1}\right|x\right)e^{-bx}\,dx \\[5pt]
&=\frac{b^{2n+1}(1+a_2-a_1)_n}{n!}
\sum_{k=0}^n\frac{(-n)_k(a_2-a_1-n)_k}{k!\,(1+a_2-a_1)_k}
\int_0^{\infty}x^{n+a_2+k}e^{-bx}\,dx \\[3pt]
&=b^{2n+1}\sum_{k=0}^n
\binom{n+a_1-a_2}{k}\binom{n+a_2-a_1}{n-k}\frac{\Gamma(n+a_2+k+1)}{b^{n+a_2+k+1}}.
\label{eq15}
\end{split}
\end{equation}
Thus we have
$$I_1=q_n(a_1,a_2,b)\frac{\Gamma(a_2+1)}{b^{a_2}},
$$
which implies (\ref{eq13}). Similarly, denoting the second integral on the right-hand side of (\ref{eq14})
by $I_2$ and integrating by parts, we have
\begin{equation}
I_2=\frac{(-1)^n b^{2n+1}}{n!^2}\int_0^{\infty}
\left(\frac{x^{a_1-a_2}}{1-x}\right)^{(n)} x^{n+a_2-a_1}(x^{n+a_1}(1-x)^{2n+1}e^{-bx})^{(n)}\,dx.
\label{eq16}
\end{equation}
Applying formula (\ref{id}) with $\alpha=a_1-a_2,$ $\beta=-1$ we get
$$
\left(\frac{x^{a_1-a_2}}{1-x}\right)^{(n)}=(a_1-a_2-n+1)_nx^{a_1-a_2-n}(1-x)^{-n-1}
{}\sb 2F\sb
 {1}\!\left(\!\left.\atop{-n,
a_1-a_2-n}{a_1-a_2-n+1}\right|x\right).
$$
Substituting the $n$-th derivative into the integral (\ref{eq16}) and integrating by parts we obtain
\begin{equation*}
\begin{split}
I_2&=\frac{b^{2n+1}(a_1-a_2-n+1)_n}{ n!^2} \\
&\times\int_0^{\infty}
\left((1-x)^{-n-1}{}\sb 2F\sb
 {1}\!\left(\!\left.\atop{-n,
a_1-a_2-n}{a_1-a_2-n+1}\right|x\right)\right)^{(n)}x^{n+a_1}(1-x)^{2n+1}e^{-bx}\,dx.
\end{split}
\end{equation*}
Using the following relation (see \cite[\S 3.4, (10)]{lu}):
$$
\frac{d^n}{dx^n}\left[(1-x)^{\alpha+\beta-\delta}{}\sb 2F\sb
 {1}\!\left(\!\left.\atop{\alpha,
\beta}{\delta}\right|x\right)\right]=\frac{(\delta-\alpha)_n(\delta-\beta)_n}{(\delta)_n}
(1-x)^{\alpha+\beta-\delta-n}{}\sb 2F\sb
 {1}\!\left(\!\left.\atop{\alpha,
\beta}{\delta+n}\right|x\right)
$$
with $\alpha=-n,$ $\beta=a_1-a_2-n,$ $\delta=a_1-a_2-n+1$ we can compute the $n$-th derivative
in the last integral
\begin{equation*}
\begin{split}
\left((1-x)^{-n-1}{}\sb 2F\sb
 {1}\!\left(\!\left.\atop{-n,
a_1-a_2-n}{a_1-a_2-n+1}\right|x\right)\right)^{(n)}&=
\frac{(a_1-a_2+1)_n\, n!}{(a_1-a_2-n+1)_n}(1-x)^{-2n-1} \\[2pt]
&\times {}\sb 2F\sb
 {1}\!\left(\!\left.\atop{-n,
a_1-a_2-n}{a_1-a_2+1}\right|x\right)
\end{split}
\end{equation*}
to get
$$
I_2=\frac{b^{2n+1}(a_1-a_2+1)_n}{n!}
\int_0^{\infty}x^{n+a_1}{}\sb 2F\sb
 {1}\!\left(\!\left.\atop{-n,
a_1-a_2-n}{a_1-a_2+1}\right|x\right)e^{-bx}\,dx.
$$
Comparing this integral with (\ref{eq15}) we obtain
$$
I_2=q_n(a_2,a_1,b) \frac{\Gamma(a_1+1)}{b^{a_1}},
$$
and the desired expansion of $R_n(a_1,a_2,b)$ is proved.

From the definition of the sequence $q_n(a_1,a_2,b)$ it follows that $q_n(a_1,a_2,b)
\in {\mathbb Q}[b].$ To estimate its denominators
we need the following property of binomial coefficients (see \cite[lemma 4.1]{chu}):
\begin{equation}
\mu_a^k \cdot\frac{(a)_k}{k!}\in {\mathbb Z} \qquad \mbox{if}\quad a\in {\mathbb Q}.
\label{eq17}
\end{equation}
Then according to (\ref{eq17}), the numbers
$$
\mu_{a_2}^{2n}\mu_{a_2-a_1}^n \binom{n+a_2-a_1}{n-k}\binom{n+a_1-a_2}{k} \frac{(a_2+1)_{n+k}}{(n+k)!}
$$
are integers for $k=0,1, \ldots, n.$ This proves the lemma. \qed

\vspace{0.2cm}

\noindent {\it Remark 1.} From Lemma \ref{l1} it follows that $R_n(a_1,a_2,b)=-R_n(a_2,a_1,b).$
Note that this symmetry property can be obtained directly from (\ref{eq12}), since we have
$Q_n^{(a_1,a_2,b)}(x)=Q_n^{(a_2,a_1,b)}(x).$ The last equality follows easily from the
Rodrigues formula (\ref{eq02}) and the fact (see \cite[Ex.]{as}) that the Rodrigues operators
corresponding to the weights $w_1$ and $w_2$ in (\ref{eq02}) commute.

Note that by  the analogous arguments, the Rivoal construction (\ref{10.5}) and therefore,
the polynomials $P_n(a,x), Q_n(a,x)$ possess a similar symmetry. Indeed, if we define
the generalized Laguerre polynomial (see \cite[\S 3.2]{as})
\begin{equation*}
\begin{split}
A_{n,a_1,a_2}(x)&=\frac{x^{-a_2}e^x}{n!^2}(x^{n+a_2-a_1}(x^{n+a_1}e^{-x})^{(n)})^{(n)}
\\ &=
\frac{(a_1+1)_n(a_2+1)_n}{n!^2}\,e^x
{}\sb 2F\sb
 {2}\left(\left.\atop{a_1+n+1,
a_2+n+1}{a_1+1, \, a_2+1}\right|-x\right)
\end{split}
\end{equation*}
and corresponding to it integrals
$$
R_{n,a_j}(x)=\int_0^{\infty}\frac{A_{n,a_1,a_2}(t)}{x-t}t^{a_j}e^{-t}\,dt
\qquad j=1,2
$$
and
$$
R_{n,a_1,a_2}(x)=x^{a_1}R_{n,a_2}(x)-x^{a_2}R_{n,a_1}(x),
$$
then we will have $A_{n,a_1,a_2}(x)=A_{n,a_2,a_1}(x),$
$R_{n,a_1,a_2}(x)=-R_{n,a_2,a_1}(x).$ This implies that the ``determinantal'' sequence
\begin{equation}
\begin{split}
{\mathcal S}_n(a_1,a_2,x)&=\begin{vmatrix}
A_{n,a_1,a_2}(x) & R_{n,a_1,a_2}(x) \\
A_{n+1,a_1,a_2}(x) & R_{n+1,a_1,a_2}(x)
\end{vmatrix} \\[3pt]
&=Q_n(a_1, a_2, x)\, x^{a_2}\Gamma(a_1+1)-P_n(a_1,a_2,x)\, x^{a_1}\Gamma(a_2+1)
\label{rem}
\end{split}
\end{equation}
has the following properties:
$$
{\mathcal S}_n(a_1,a_2,x)=-{\mathcal S}_n(a_2,a_1,x) \qquad\mbox{and}\qquad
P_{n}(a_1, a_2, x)=Q_n(a_2, a_1,x).
$$
Taking $a_1=a,$ $a_2=0$ in (\ref{rem}) we get (\ref{10.5}), which yields
$$
Q_n(a,x)=Q_n(a,0,x), \qquad P_n(a,x)=Q_n(0,a,x).
$$
It would be of interest to find a compact hypergeometric form for the polynomial
$Q_n(a_1,a_2,x),$ if it exists (maybe as a product of hypergeometric series),
for better understanding of  hypergeometric origin of such sequences.

\vspace{0.2cm}

The next lemma provides an alternative representation of $R_n(a_1,a_2,b)$ as a multiple integral.
\begin{lemma} \label{l2}
Let $a_1, a_2, b\in {\mathbb R},$ $a_2-a_1\not\in {\mathbb Z},$ $a_1, a_2>-1,$ $b>0.$ Then one has
\begin{equation*}
\begin{split}
R_n(a_1,a_2,b)&=(a_2-a_1)
\frac{(1+a_2-a_1)_{n}(1+a_1-a_2)_n \,b^{2n+1}}{n!^2} \\[2pt]
&\times
\int_0^{\infty}\!\!\!\int_0^{\infty}\!\!
\frac{x^{n+a_1}y^n(x-1)^{2n+1}e^{-bx}}{(xy+1)^{n+1+a_1-a_2}(y+1)^{n+1+a_2-a_1}}\,
dxdy.
\end{split}
\end{equation*}
\end{lemma}
{\bf Proof.} Substituting (\ref{eq02}) into (\ref{eq12}) we obtain
$$
R_n(a_1,a_2,b)=\frac{b^{2n+1}}{n!^2}\int_0^{\infty}
\frac{1-x^{a_1-a_2}}{1-x}\bigl(x^{n+a_2-a_1}(x^{n+a_1}(1-x)^{2n+1}e^{-bx})^{(n)}\bigr)^{(n)}\,dx.
$$
Applying the following representation introduced by Rivoal (see \cite[proof of Lemma 4]{ri1}):
$$
\frac{1}{a_1-a_2}\cdot\frac{1-x^{a_1-a_2}}{1-x}=\int_0^{\infty}
\frac{du}{(1+ux)^{1+a_2-a_1}(1+u)^{1+a_1-a_2}}
$$
we get
$$
R_n(a_1,a_2,b)=\frac{(a_1-a_2) b^{2n+1}}{n!^2}
\int_0^{\infty}\!\!\int_0^{\infty}\frac{\bigl(x^{n+a_2-a_1}(x^{n+a_1}(1-x)^{2n+1}e^{-bx})^{(n)}
\bigr)^{(n)}}{(1+ux)^{1+a_2-a_1}(1+u)^{1+a_1-a_2}}\,dudx.
$$
Now using Fubini's theorem and $n$-integrations by parts with respect to $x$ we obtain
$$
R_n(a_1,a_2,b)=\frac{(a_2-a_1)_{n+1} b^{2n+1}}{n!^2}
\int_0^{\infty}\!\!\int_0^{\infty}\frac{u^nx^{n+a_2-a_1}(x^{n+a_1}(x-1)^{2n+1}e^{-bx})^{(n)}}%
{(1+ux)^{n+1+a_2-a_1}(1+u)^{1+a_1-a_2}}dudx.
$$
Making the change of  variable $u=1/(xy)$ or $y=1/(xu)$ in the last integral we get
$$R_n(a_1,a_2,b)=\frac{(a_2-a_1)_{n+1} b^{2n+1}}{n!^2}
\int_0^{\infty}\!\!\int_0^{\infty}\frac{(x^{n+a_1}(x-1)^{2n+1}e^{-bx})^{(n)}}%
{(y+1)^{n+1+a_2-a_1}(xy+1)^{1+a_1-a_2}}\,dxdy.
$$
Finally, integrating by parts with respect to $x$ we obtain the desired representation
and the lemma is proved. \qed

Our next concern will be the asymptotic behavior of the  sequences $q_n(a_1,a_2,b)$ and
$R_n(a_1,a_2,b).$
\begin{lemma} \label{l3}
Let $a_1, a_2, b\in {\mathbb R},$ $a_2-a_1\not\in {\mathbb Z},$ $a_1, a_2>-1,$ $b>0.$
Then the following asymptotic formulas are valid:
\begin{equation*}
\begin{split}
R_n(a_1,a_2,b)&=(2n)! \frac{e^{-\sqrt{2bn}}}{n^{1/4-(a_1+a_2)/2}}\bigl(
2\sin(\pi(a_2-a_1))\frac{\Gamma(a_2+1)}{b^{a_2}}c_2+O(n^{-1/2})\bigr), \\[2pt]
q_n(a_1,a_2,b)&=(2n)! \frac{e^{\sqrt{2bn}}}{n^{1/4-(a_1+a_2)/2}}\bigl(c_2+O(n^{-1/2})\bigr)
\quad \mbox{as} \quad n\to\infty,
\end{split}
\end{equation*}

\vspace{0.2cm}

\noindent where the constant $c_2$ is defined in Theorem {\rm 1}.
\end{lemma}
{\bf Proof.} To obtain the asymptotic behavior of the above sequences,  we proceed analogously to the proof
of Lemma 5 in \cite{he}.
Substituting (\ref{eq02}) into (\ref{eq13}) and integrating by parts  we have
$$
q_n(a_1,a_2,b)=\frac{(-1)^n}{n!}\frac{b^{2n+a_2+1}}{\Gamma(a_2+1)}\int_0^{\infty}
x^{n+a_1}(x-1)^{2n+1}e^{-bx}\left(\frac{x^{n+a_2-a_1}}{(x-1)^{n+1}}\right)^{(n)}\,dx.
$$
Further, repeating the arguments such as those used in \cite[p.56-59]{aly} we get
$$
q_n(a_1,a_2,b)=\frac{b^{2n+a_2+1}}{\Gamma(a_2+1)}\frac{n^{2n+(a_1+a_2)/2+3/4}}{2\sqrt{\pi}}
\int_0^{\infty}e^{n\varphi_1(y)}\psi_1(y)\,dy\cdot (1+O(n^{-1/2})),
$$
where
$$
\varphi_1(y)=2\ln(\sqrt{y}+1/\sqrt{n})+\ln y- by, \qquad
\psi_1(y)=(\sqrt{y}+1/\sqrt{n})y^{(a_1+a_2)/2-1/4}.
$$
The last integral coincides with the corresponding integral from \cite[Lemma 5]{he} and likewise, by the Laplace method, we obtain
\begin{equation*}
\begin{split}
q_n(a_1,a_2,b)&=\frac{(2n)^{2n+(a_1+a_2)/2+1/4}}{\Gamma(a_2+1)b^{1/4+(a_1-a_2)/2}}e^{-2n+\sqrt{2bn}-3b/8}
\bigl(1+O(n^{-1/2})\bigr) \\[2pt]
&=(2n)! n^{(a_1+a_2)/2-1/4}e^{\sqrt{2bn}}\bigl(c_2+O(n^{-1/2})\bigr),
\end{split}
\end{equation*}
where the constant $c_2=c_2(a_1,a_2,b)$ is defined in Theorem \ref{t1}.

Making the change of  variables $x=nu, y=v/\sqrt{nu}$ in the double integral of Lemma~\ref{l2}
we get
$$
R_n
=\frac{\Gamma(n+1+a_2-a_1)\Gamma(n+1+a_1-a_2)}{n!^2\Gamma(a_2-a_1)\Gamma(1+a_1-a_2)}
b^{2n+1}n^{3n+2+(a_1+a_2)/2}\!\!\int_0^{\infty}\!\!\!g(u,v) f^n(u,v)\,dudv,
$$
where
$$
f(u,v)=\frac{uv(u-1/n)^2e^{-bu}}{(v+\sqrt{nu})(v\sqrt{nu}+1)},\qquad
g(u,v)=\frac{u^{(a_1+a_2)/2}(u-1/n)}{(v+\sqrt{nu})^{1+a_2-a_1}
(v\sqrt{nu}+1)^{1+a_1-a_2}}.
$$
Using the reflection formula
$$
\Gamma(a_2-a_1)\cdot \Gamma(1+a_1-a_2)=\frac{\pi}{\sin(\pi(a_2-a_1))}
$$
and Stirling's formula for the gamma function
$$
\Gamma(n+1+a_2-a_1)=n^{a_2-a_1} n! \bigl(1+O(n^{-1/2})\bigr) \quad\mbox{as}\quad n\to\infty
$$
we get
$$
R_n=\frac{\sin(\pi(a_2-a_1))}{\pi}b^{2n+1}n^{3n+2+(a_1+a_2)/2}\int_0^{\infty}g(u,v)f^n(u,v)\,dudv
\cdot\bigl(1+O(n^{-1})\bigr).
$$
Now applying the Laplace method for multiple integrals by the same way as in the proof of Lemma 5
\cite{he} with the only difference that in our case
$$
g(u_0,v_0)=\frac{1}{n}\left(\frac{2}{b}\right)^{\frac{a_1+a_2}{2}}\bigl(1+O(n^{-1/2})\bigr),
$$
we obtain
\begin{equation*}
\begin{split}
R_n&=2\sin(\pi(a_2-a_1))\frac{(2n)^{2n+1/4+(a_1+a_2)/2}}{b^{1/4+(a_1+a_2)/2}}e^{-2n-\sqrt{2bn}-3b/8}
\bigl(1+O(n^{-1/2})\bigr) \\
&=(2n)!e^{-\sqrt{2bn}}n^{(a_1+a_2)/2-1/4}
\bigl(2\sin(\pi(a_2-a_1))\frac{\Gamma(a_2+1)}{b^{a_2}} c_2+O(n^{-1/2})\bigr)
\end{split}
\end{equation*}
and the lemma is proved. \qed

Now Theorem \ref{t1} follows easily from Lemmas \ref{l1}--\ref{l3}.

\section{Proof of Theorem \ref{t2}.}

It is readily seen that the construction (\ref{pn}) is a limiting case of Theorem \ref{t1}.
Indeed, comparing Lemma \ref{l2} with Lemmas 3, 4 from \cite{he} we get
\begin{equation}
\lim_{a_2\to a_1}\frac{b^{a_1}}{\Gamma(a_1+1)}\frac{R_n(a_1,a_2,b)}{a_2-a_1}=
p_n(a_1,b)-(\ln b-\psi(a_1+1))q_n(a_1,b).
\label{eq18}
\end{equation}
On the other hand, by Lemma \ref{l1}, we have
\begin{equation*}
\begin{split}
&\lim_{a_2\to a_1}\frac{b^{a_1}}{\Gamma(a_1+1)}\frac{R_n(a_1,a_2,b)}{a_2-a_1}=
\lim_{a_2\to a_1}\frac{1}{a_2-a_1}\biggl(q_n(a_1,a_2,b)\frac{\Gamma(a_2+1)/b^{a_2}}{\Gamma(a_1+1)/b^{a_1}}
-q_n(a_2,a_1,b)\biggr)\\[3pt]
&=\lim_{a_2\to a_1}\frac{q_n(a_1,a_2,b)}{\Gamma(a_1+1)}\cdot \frac{\frac{\Gamma(a_2+1)}{b^{a_2-a_1}}-\Gamma(a_1+1)}%
{a_2-a_1}+\lim_{a_2\to a_1} \frac{q_n(a_1,a_2,b)-q_n(a_2,a_1,b)}{a_2-a_1}\\[3pt]
&=\frac{q_n(a_1,a_1,b)}{\Gamma(a_1+1)}\cdot\frac{d}{da_2}\left.\biggl(\frac{\Gamma(a_2+1)}{b^{a_2-a_1}}\biggr)
\right|_{a_2=a_1}+\lim_{a_2\to a_1} \frac{q_n(a_1,a_2,b)-q_n(a_2,a_1,b)}{a_2-a_1}.
\end{split}
\end{equation*}
Taking into account that
$$
\frac{d}{da_2}\biggl(\frac{\Gamma(a_2+1)}{b^{a_2-a_1}}\biggr)=\frac{\Gamma(a_2+1)}{b^{a_2-a_1}}(\psi(a_2+1)
-\ln b),
$$
we obtain
\begin{equation}
\lim_{a_2\to a_1}\!\frac{b^{a_1}}{\Gamma(a_1+1)}\frac{R_n(a_1,a_2,b)}{a_2-a_1}=
q_n(a_1,b)(\psi(a_1+1)-\ln b)+\lim_{a_2\to a_1}\!\!\frac{q_n(a_1,a_2,b)-q_n(a_2,a_1,b)}{a_2-a_1}.
\label{eq19}
\end{equation}
Comparing the right-hand sides of (\ref{eq18}) and  (\ref{eq19}) we conclude that
$$
p_n(a_1,b)=\lim_{a_2\to a_1} \frac{q_n(a_1,a_2,b)-q_n(a_2,a_1,b)}{a_2-a_1}.
$$
To evaluate the above limit, we need the following expansions:
\begin{align*}
\binom{n+a_1-a_2}{k}&=\binom{n}{k}-\binom{n}{k}(H_n-H_{n-k})(a_2-a_1)+
O\bigl((a_2-a_1)^2\bigr), \\[2pt]
\binom{n+a_2-a_1}{n-k}&=\binom{n}{k}+\binom{n}{k}(H_n-H_k)(a_2-a_1)+O\bigl((a_2-a_1)^2\bigr),\\[4pt]
(a_2+1)_{n+k}&=(a_1+1)_{n+k}+(a_1+1)_{n+k}H_{n+k}(a_1)(a_2-a_1)+O\bigl((a_2-a_1)^2\bigr).
\end{align*}
Then a straightforward verification shows that
$$
p_n(a_1,b)=\sum_{k=0}^n b^{n-k}\binom{n}{k}^2(a_1+1)_{n+k}(H_{n+k}(a_1)+2H_{n-k}-2H_k),
$$
and this is precisely the assertion of Theorem \ref{t2}. \qed

Setting $a_1=0,$ $b=1$ yields Corollary \ref{c1}.

\vspace{0.2cm}

\noindent {\it Remark 2.} Note that an analogous conclusion can be drawn for the polynomials $P_n(x), Q_n(x)$
from (\ref{10.1}). Comparing \cite[Lemma 1]{ri} with \cite[Lemma 4]{ri1} gives
$$
{\mathcal S}_n(x)=\lim_{a\to 0}\,\frac{-1}{a\,x^a}\cdot{\mathcal S}_n(a,x).
$$
Then analysis similar to that in the proof of Theorem \ref{t2} shows that
$Q_n(x)=Q_n(0,x)=Q_n(0,0,x)$ and
$$
P_n(x)=\lim_{a\to 0} \frac{Q_n(a,x)-P_n(a,x)}{a}=
\lim_{a\to 0}\frac{Q_n(a,0,x)-Q_n(0,a,x)}{a},
$$
where the polynomial $Q_n(a_1,a_2,x)$ is defined in Remark 1.

\end{document}